\documentclass[12pt]{article}

\usepackage{graphicx,psfrag}
\usepackage{psfrag,eepic,epsfig}
\usepackage{sectsty}
\usepackage{setspace}
\usepackage{color}

\usepackage{setspace}
\usepackage{amsmath, epsfig, cite, ifpdf}
\usepackage{amssymb}
\usepackage{amsfonts}
\usepackage{latexsym}
\usepackage{graphicx,psfrag,eepic,rawfonts}
\usepackage{enumerate}
\usepackage{indentfirst}
\input{epsf}

\newtheorem{thm}{Theorem}[section]
\newtheorem{prop}[thm]{Proposition}

\newtheorem{defi}[thm]{Definition}
\newtheorem{lem}[thm]{Lemma}
\newtheorem{conj}[thm]{Conjecture}

\newcommand{\qed}{{\hfill\rule{4pt}{7pt}}}
\def\pf{\noindent {\it Proof.} }

\numberwithin{equation}{section}

\makeatletter \@addtoreset{equation}{section} \makeatother

\title {\bf Solution to a conjecture on the maximum skew-spectral
radius of odd-cycle graphs\footnote{Supported by NSFC No.11371205,
the 973 program No.2013CB834204, and PCSIRT.}}

\author{
{\small Xiaolin Chen, Xueliang Li}\\
{\small Center for Combinatorics and LPMC-TJKLC}\\
{\small Nankai University, Tianjin 300071, P.R. China}\\
{\small E-mail: chxlnk@163.com; lxl@nankai.edu.cn}\\
{\small Huishu Lian}\\
{\small College of Science}\\
{\small China University of Mining and Technology, Xuzhou 221116, P.R. China}\\
{\small E-mail:lhs6803@126.com}   }
\date{}
\begin{document}

\maketitle

\begin{abstract}
Let $G$ be a simple graph with no even cycle, called an odd-cycle
graph. Cavers et al. [Cavers et al. Skew-adjacency matrices of
graphs, Linear Algebra Appl. 436(2012), 4512--1829] showed that the
spectral radius of $G^\sigma$ is the same for every orientation
$\sigma$ of $G$, and equals the maximum matching root of $G$. They
proposed a conjecture that the graphs which attain the maximum skew
spectral radius among the odd-cycle graphs $G$ of order $n$ are
isomorphic to the odd-cycle graph with one vertex degree $n-1$ and
size $m=\lfloor 3(n-1)/2\rfloor$. This paper, by using the Kelmans
transformation, gives a proof of the conjecture. Moreover, sharp
upper bounds of the maximum matching roots of the odd-cycle graphs
with given order $n$ and size $m$ are given and extremal graphs are
characterized.
\\


\end{abstract}

\section{Introduction}

Let $G$ be a finite simple graph with vertex set
$V=\{v_1,\ldots,v_n\}$ and edge set $E(G)$. For more notation and
terminology that will be used in the sequel, we refer to
\cite{BM,BH}. Given an orientation $\sigma$ of $G$ which makes any
edge $uv$ to be an arc $u\rightarrow v$ or $v\rightarrow u$, we get
an oriented graph $G^\sigma$. The skew-adjacency matrix is the
square matrix $S(G^\sigma)=[s_{ij}]$ of order $n$, where $s_{ij}=1$
and $s_{ji}=-1$ if $u\rightarrow v$, $s_{ij}=s_{ji}=0$ otherwise. By
the definition, $S(G^\sigma)$ is a skew symmetric matrix. Let the
skew-characteristic polynomial be
$\phi_{s}(G^\sigma,x)=det(xI_n-S(G^\sigma))$, where $I_n$ is the
unit matrix with order $n$. Thus the eigenvalues of $S(G^\sigma)$
are all pure imaginary numbers, which form the skew-spectrum of
$G^\sigma$. The spectral radius $\rho(G^\sigma)$ of an oriented
graph $G^\sigma$ is the maximum norm of its all eigenvalues. In
recent years, there are many results about the skew-spectrum and the
skew-spectral radius, see \cite{ABS,CCF,CLL,SS}.

For a given graph $G$, different orientations may lead to different
skew-spectra, and thus different skew-spectral radius. In
\cite{CCF}, Cavers et al. introduced the maximum skew-spectral
radius $\rho_s(G)$ which is defined as $\rho_s(G)=
\max\{\rho(G^\sigma):\sigma$ is an orientation of $G$\}. They
studied a special class of graphs, called the odd-cycle graphs. A
graph is called an odd-cycle graph if it have no even cycles. In
\cite{CCF}, Cavers et al. showed that when the graph $G$ is an
odd-cycle graph, the skew-spectra are independent of the orientation
$\sigma$. In fact, the skew-characteristic polynomial
$\phi_{s}(G^\sigma,x)$ of $G^\sigma$ can be determined by the
matching polynomial $m(G,x)$ of the graph $G$. Anuradha and
Balakrishnan also got some skew-spectral properties of the odd-cycle
graphs in \cite{ABS}.

Now, we recall the definition of the matching polynomial \cite{GG}.
\begin{defi}
Let $m_r(G)$ denote the number of $r$ independent edges in a graph
$G$. Define the matching polynomial of $G$ as
\begin{equation*}
m(G,x)=\sum_{k=0}(-1)^{k}m_k(G)x^{n-2k}.
\end{equation*}
\end{defi}

Note that many results about the roots of the matching polynomial
have been obtained; see \cite{FA,GOD,GG,HL}, such as, the matching
roots are real-valued, and symmetric. Denote by $t(G)$ the maximum
matching root of $G$. For the applications of the matching
polynomial in chemistry and statistical physics, we refer the reader
to \cite{HL,GMT}.

Given an odd-cycle graph $G$, the following lemma in \cite{CCF}
implies the relationship between the skew-characteristic polynomial
of an oriented graph $G^\sigma$ and the matching polynomial of the
graph $G$.
\begin{lem}{(\cite{CCF})}\label{Skewspectral}
Let $G$ be a graph of order $n$. Then $G$ is an odd-cycle graph if
and only if $\phi_s(G^\sigma,x)=(-i)^nm(G,ix)$ for all orientations
$\sigma$ of $G$.
\end{lem}

From the above lemma, we note that when $G$ is an odd-cycle graph,
the skew-spectral radius $\rho(G^\sigma)$ under any orientation
$\sigma$ equals the maximum skew-spectral radius $\rho_s(G)$.
Moreover, the maximum skew-spectral radius of the graph $G$ equals
the maximum matching root of the graph $G$, i.e. $\rho_s(G)=t(G)$.

In \cite{CCF}, Cavers et al. studied the upper bound of the maximum
skew-spectral radius among the odd-cycle graphs with order $n$, and
proposed a conjecture. Let $H_n$ be the odd-cycle graph of order $n$
with one vertex degree $n-1$ and size $m=\lfloor 3(n-1)/2\rfloor$.
It is easy to check that $H_n$ is unique (up to isomorphism).
Moreover, we can find that $H_n$ is the maximal odd-cycle graph of
order $n$, since the size $m$ of an odd-cycle graph of order $n$ is
no more than $\lfloor 3(n-1)/2\rfloor$ (\cite{ABS}, Theorem 3.2).
\begin{conj}\label{conjecture} (\cite{CCF})
If $G$ is an odd-cycle graph of order $n$, then $\rho_s(G)\leq
\rho_s(H_n)$, and equality holds if and only if $G\cong H_n$.
\end{conj}

Based on Lemma \ref{Skewspectral}, we can find that Conjecture
\ref{conjecture} can be equivalently rewritten as follows.
\begin{conj}\label{conjecture2}
If $G$ is an odd-cycle graph of order $n$, then $t(G)\leq t(H_n)$,
and equality holds if and only if $G\cong H_n$.
\end{conj}

For convenience, in the remainder of this paper, we just study the
maximum matching root $t(G)$ of an odd-cycle graph. We will prove
Conjecture \ref{conjecture2}. In fact, we get an even stronger
result. We characterize the extremal graphs with maximum $t(G)$
among the odd-cycle graphs with order $n$ and size $m$.
\begin{thm}\label{THM3}
Let $G$ be an odd-cycle graph with order $n$ and size $m$, where
$1\leq m\leq \lfloor \frac{3(n-1)}{2}\rfloor$. Let $F(n,m)$ be an
odd-cycle graph with one vertex degree $n-1$ and size $m$, which is
unique up to isomorphism.
\begin{enumerate}[(1)]
\item If $m=1$, then $t(G)=1$.
\item If $m=2$, then $t(G)\leq \sqrt{2}$, and equality holds if and only if $G$ is the
disjoint union of a star $K_{1,2}$ and $n-2$ isolated vertices.
\item If $m=3$, then $t(G)\leq \sqrt{3}$, and equality holds if and only if $G$ is the
disjoint union of a triangle and $n-3$ isolated vertices, or the
disjoint union of a star graph $K_{1,3}$ and $n-4$ isolated
vertices.
\item If $4\leq m\leq n-2$, then $t(G)\leq t(F(m+1,m))$,
and equality holds if and only if $G$ is the disjoint union of the
graph $F(m+1,m)$ and $n-m-1$ isolated vertices.
\item If $n-1\leq m\leq \lfloor \frac{3(n-1)}{2}\rfloor$ and $m\geq 4$, then
$t(G)\leq t(F(n,m))$, and equality holds if and only if $G\cong
F(n,m)$.
\end{enumerate}
\end{thm}

We organize the rest of this paper as follows. In Section
\ref{sec2}, we study the Kelmans transformation acting on odd-cycle
graphs. In Section \ref{sec3}, we find the monotone property of the
maximum matching roots after an acting of the Kelmans
transformation. In Section \ref{sec4}, we give our proofs of
Conjecture \ref{conjecture2} and Theorem \ref{THM3}.

\section{The Kelmans transformation on odd-cycle graphs}\label{sec2}

In \cite{CCF}, Cavers et al. suggested that the standard techniques
called edge-switching will be useful to prove Conjecture
\ref{conjecture}. In this paper, we mainly use a graph
transformation called the Kelmans transformation.

\begin{figure}[h,t,b,p]
\begin{center}
\scalebox{0.7}[0.7]{\includegraphics{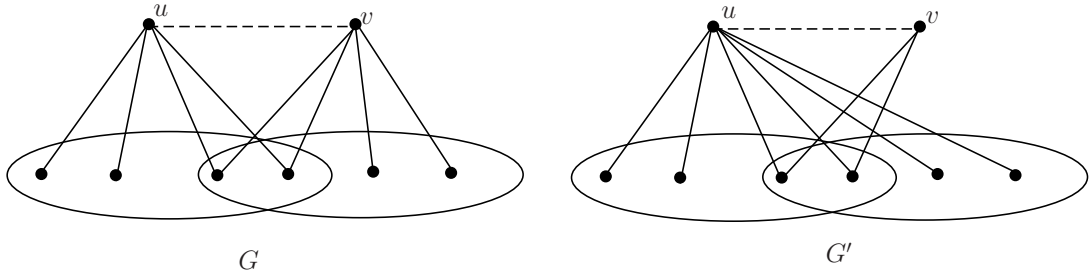}}
\end{center}
\caption{The Kelmans transformation.}\label{Fig1}
\end{figure}

We first recall some results about the Kelmans transformation.
In\cite{KEL}, Kelmans introduced a graph transformation for the
extremal problems related with the synthesis of reliable networks.
In \cite{SSS}, Satyanarayana et al. also introduced a reliability
improving graph transformation, named the "swing surgery", which can
be thought as the inverse of the Kelmans transformation. Brown et
al. \cite{BCD} called the Kelmans transformation as shift operation,
and studied many applications in network reliability. Csikv\'{a}ri
in \cite{PC1} made a breakthrough in solving a conjecture of
Nikiforov by the Kelmans transformation. Moreover, Csikv\'{a}ri in
\cite{PC} found several applications of the Kelmans transformation
in the extremal problems on graph polynomials, such as, matching
polynomial, independent polynomial, chromatic polynomial, Laplacian
polynomial. For more interesting details, see \cite{PC}. Next, we
give the definition of the Kelmans transformation.
\begin{defi}
Let $u$ and $v$ be two vertices of a graph $G$. Let $N(v)$ be the
neighbor set of $v$ in $G$. For any $w\in N(v)$, if $w$ is not $u$
and $u$ and $w$ are nonadjacent, then we delete the edge $vw$ and
add the edge $uw$. After the above operations, we get a new graph
$G'$. Then $G'$ is a graph obtained from $G$ by the Kelmans
transformation between $u$ and $v$; see Fig. \ref{Fig1}.
\end{defi}

For convenience, Csikv\'{a}ri \cite{PC} called $u$ and $v$ the
beneficiary and the co-beneficiary of the transformation,
respectively. For brevity, we denote by $KT(G,u,v)$ the Kelmans
transformation acting on a graph $G$ with the beneficiary $u$ and
the co-beneficiary $v$. It is easy to check that the Kelmans
transformation does not change the number of edges.

Let $G$ be a connected odd-cycle graph with order $n$ and size $m$.
Then $n-1\leq m\leq n+\lfloor(n-1)/2\rfloor-1$ (\cite{ABS}, Theorem
3.2). Let $F(n,m)$ be an odd-cycle graph with one vertex degree
$n-1$ and size $m$. It is easy to check that $F(n,m)$ is unique (up
to isomorphism). If $m=n+\lfloor(n-1)/2\rfloor-1$, then $F(n,m)$ is
the same as $H_n$.
\begin{thm}\label{KTOG}
Let $G$ be a connected odd-cycle graph with order $n$ and size $m$.
Then $F(n,m)$ can be obtained from $G$ by a number of Kelmans
transformations.
\end{thm}
\pf Let $G$ be a connected odd-cycle graph with order $n$ and size
$m$, which implies that $n-1\leq m\leq n+\lfloor(n-1)/2\rfloor-1$.
Then $G$ must be a cactus, that is, every block of $G$ is either a
cycle or an edge (\cite{BM}, EX. 3.2.3). It follows that any two
cycles of $G$ are edge-disjoint. Let $C(G)$ be the set of odd cycles
with order larger than $3$. Next, we will prove the theorem by three
steps.

\begin{figure}[h,t,b,p]
\begin{center}
\scalebox{0.7}[0.7]{\includegraphics{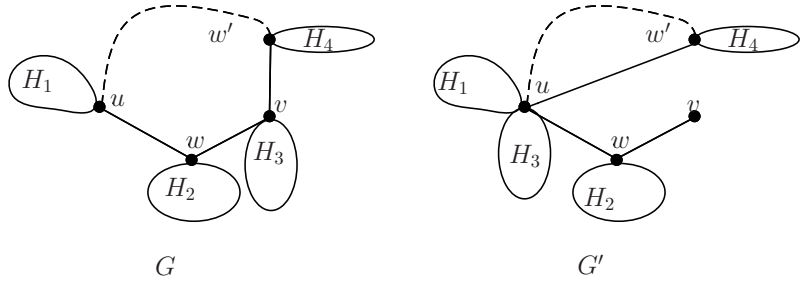}}
\end{center}
\caption{$KT(G,u,v)$ in step 1.}\label{Fig2}
\end{figure}

\textbf{Step 1.} If $C(G)$ is nonempty, assume that
$C(G)=\{C_1,\ldots,C_t\}$. The number of the edges of $C(G)$ is
$|E(C(G))|=\sum_{i=1}^t|C_t|\leq m$. Let $C\in C(G)$ and $u,v$ be
two vertices of $C$ whose distance in $C$ is $2$. Then there exists
a path $uwvw'$ in $C$. By the Kelmans transformation $KT(G,u,v)$, we
get a graph $G'$. It is easy to find that the degree of the vertex
$v$ in $G'$ is $1$; see Fig. \ref{Fig2}.

We first show that $G'$ is connected. For every vertex $u'$ in $G$,
there is a path $u'Pu$ in $G$ connected $u'$ and $u$. If the vertex
$v$ is in the path $P$, then the path $P$ can be decomposed as
$P_1vP_2$. By the definition of the Kelmans transformation, we can
check that $u'P_1u$ is in the graph $G'$ and $u'$ is connected to
$u$ in $G'$. If $v$ is not in the path $P$, then $u'$ and $u$ are
connected by the path $P$ in $G'$.

Secondly, we claim that $G'$ is still an odd-cycle graph. Suppose
that there exists an even cycle $C'$ in $G'$. Then $u$ must be in
$C'$. Let $C'=uh_1\ldots h_su$. If $uh_1$ and $uh_s$ both belong to
$G$, then $C'$ is also an even cycle of $G$, a contradiction. If
neither $uh_1$ nor $uh_s$ belongs to $G$, then $vh_1\ldots h_sv$ is
an even cycle of $G$, again a contradiction. Therefore, without loss
of generality, suppose that $uh_1$ is in $G$ and $uh_s$ is not in
$G$, which implies that $h_1\neq h_s$ and $vh_s$ is an edge of $G$.
We then obtain a path $P=uh_1\cdots h_sv$ with even length in $G$.
Since $G$ is a cactus and $u,v$ are both in the cycle $C$, we can
deduce that $P$ is in the cycle $C$. It follows that $P=uwv$,
contradicting with $h_1\neq h_s$.

Finally, we consider the cycle set $C(G')$. By analyzing the
corresponding relationship between $C(G)$ and $C(G')$, it can be
verified that $|C(G')|=|C(G)|\,\,\mbox{or}\,\,|C(G)|-1$ and
$|E(C(G'))|\leq|E(C(G))|-2$. Thus, by less than the number
$\frac{m}{2}$ of Kelmans transformations as above, we can get a
connected odd-cycle graph $G'$ whose cycles are all triangles. Then,
we turn to Step 3.

\textbf{Step 2.} If $C(G)$ is empty, then $G$ contains no cycles or
each cycle of $G$ is a triangle. Let $G':=G$.

\begin{figure}[h,t,b,p]
\begin{center}
\scalebox{0.7}[0.7]{\includegraphics{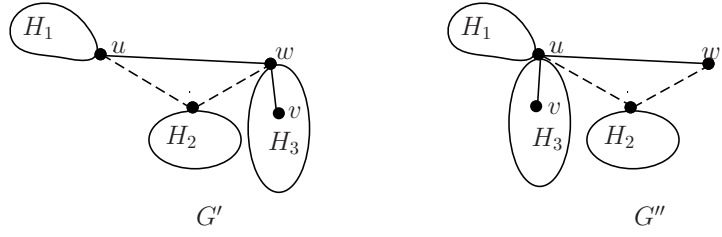}}
\end{center}
\caption{$KT(G',u,w)$ in step 3.}\label{Fig3}
\end{figure}

\textbf{Step 3.} Assume that a vertex $u$ of $G'$ has the maximum
degree in $G'$ and its degree is denoted by $d_{G'}(u)$. If
$d_{G'}(u)$ is less than $n-1$, then there exists a vertex $v$ such
that the distance of $u,v$ in $G'$ is $2$. Suppose that $uwv$ is a
$2$-path in $G'$. We apply the Kelmans transformation $KT(G',u,w)$
and get a graph $G''$. By the similar analysis in the first step, we
find that $G''$ is a connected odd-cycle graph whose cycles are all
triangles. Moreover, we get $d_{G'}(u)\geq d_{G}(u)+1$. Then by at
most $n-1$ kelmans transformations as above, we get a graph $G''$
with the maximum degree $n-1$, order $n$ and size $m$. Moreover,
$G''$ satisfies that all cycles are triangles. It is easy to see
that the graph $G''$ is isomorphic to the graph $F(n,m)$; see Fig.
\ref{Fig3}.

The proof is this complete. \qed

\section{Maximum matching root and Kelmans transformation}\label{sec3}

In this section, we will show that the maximum matching roots of the
graph strictly increase after the Kelmans transformation. Before
proceeding, we first review some useful lemmas of the matching
polynomial.
\begin{lem}(\cite{GG})\label{deledge}
Let edge $e=uv$. Then
\begin{equation}
m(G,x)=m(G-e,x)-m(G-u-v,x)
\end{equation}
\end{lem}

Let $G$ and $H$ be two disjoint graphs. Then the graph $G\bigcup H$
denotes the union of $G$ and $H$.
\begin{lem}(\cite{GG})\label{union}
Let $G_1,\ldots,G_k$ be $k$ disjoint graphs. Then
$$m(\bigcup_{i=1}^kG_i,x)=\prod_{i=1}^km(G_i,x).$$
\end{lem}

In \cite{GUT}, Gutman showed some parallel results for the roots of
the matching polynomial and the spectra of the characteristic
polynomial.
\begin{lem}(\cite{GUT})\label{sub1}
If $H$ is a subgraph of $G$, then $t(G)\geq t(H)$. If $G$ is
connected and $H$ is a proper subgraph of $G$, then $t(G)>t(H)$.
\end{lem}

Now, we give some propositions of the maximum matching root under
the Kelmans transformation.
\begin{defi}
Define $G_1\succeq G_2$ if for all $x\geq t(G_1)$ we have
$m(G_2,x)\geq m(G_1,x)$. Define $G_1\succ G_2$ if for all $x\geq
t(G_1)$ we have $m(G_2,x)>m(G_1,x)$.
\end{defi}

It is easy to check that $G_1\succ G_2$ can deduce $G_1\succeq G_2$.
Note that Csikv\'{a}ri \cite{PC} introduced the definition $\succ$
which has the same meaning as our definition $\succeq$, and they
obtained many good results which inspire us. In the following, we
obtain some results about $\succ$ and $\succeq$.

Before proceeding, we list some useful results about $\succeq$ which
were proved in \cite{PC}.
\begin{prop}(\cite{PC})\label{Ptrans1}
The relation $\succeq$ is transitive. If $G_1\succeq G_2$, then
$t(G_1)\geq t(G_2)$.
\end{prop}
\begin{prop}(\cite{PC})\label{Psub}
If $H$ is a spanning subgraph of $G$, then $G\succeq H$.
\end{prop}
\begin{thm}(\cite{PC})\label{PCT}
Assume that $G'$ is a graph obtained from $G$ by a number Kelmans
transformations, then $G'\succeq G$; in particular, $t(G')\geq
t(G)$.
\end{thm}

Next, we give our results.
\begin{prop}\label{trans1}
The relation $\succ$ is transitive. If $G_1\succ G_2$, then $t(G_1)>t(G_2)$.
\end{prop}
\pf Suppose $G_1\succ G_2$. Since $t(G_1)$ is the maximum root of
the matching polynomial, it implies that for $x\geq t(G_1)$, we have
$m(G_2,x)>m(G_1,x)\geq 0$. Due to the fact that the leading
coefficient of the matching polynomial is 1, we get $t(G_1)>t(G_2)$.
If $G_1\succ G_2\succ G_3$, then $m(G_3,x)>m(G_2,x)>m(G_1,x)$ for
$x\geq max\{t(G_1),t(G_2))\}=t(G_1)$. It follows that $G_1\succ
G_3$.

\begin{prop}\label{trans2}
If $G_1\succ G_2$ and $G_2\succeq G_3$, then $G_1\succ G_3$.
If $G_1\succeq G_2$ and $G_2\succ G_3$, then $G_1\succ G_3$.
\end{prop}
\pf We only prove the first part, the second part can be proved
similarly. For $x\geq t(G_1)$, we have $m(G_2,x)>m(G_1,x)$. For
$x\geq t(G_2)$, we have $m(G_3,x)\geq m(G_2,x)$. By Propositions
\ref{Ptrans1} and \ref{trans1}, it follows that $t(G_1)>t(G_2)\geq
t(G_3)$. For $x\geq t(G_1)$, we get $m(G_3,x)>m(G_1,x)$, i.e.
$G_1\succ G_3$. \qed

\begin{prop}\label{subroot}
If $G$ is connected and $H$ is a proper spanning subgraph of $G$, then $G\succ H$.
\end{prop}
\pf Since $H$ is a proper spanning subgraph of $G$, suppose $H$ is
obtained by deleting $l$ edges $\{e_1,\ldots,e_l\}$ from $G$. Then,
$m(H,x)=m(G-e_1-\cdots-e_l,x)$. Suppose $e_1$ is $uv$. By Lemma
\ref{deledge}, we have
\begin{equation}
m(G,x)=m(G-e_1,x)-m(G-u-v,x),
\end{equation}
Since $G$ is connected, by Lemma \ref{sub1}, we have $t(G-u-v)<t(G)$
and $t(G-e_1)<t(G)$. Since the leading coefficient of the matching
polynomial is 1, it follows that for $x\geq t(G)$, we get
\begin{equation}
m(G-e_1,x)-m(G,x)=m(G-u-v,x)>0.
\end{equation}
Thus we get $G\succ G-e_1$. Since $H$ is a spanning subgraph of
$G-e_1$, by Proposition \ref{Psub}, we have $G-e_1\succeq H$. By
Proposition \ref{trans2}, we deduce that $G\succ H$. \qed

\begin{prop}\label{Proper}
If $H$ is a proper spanning subgraph of $G$, then for $x>t(G)$, we
have $m(H,x)>m(G,x)$.
\end{prop}
\pf If $G$ is a connected graph, then by Proposition \ref{subroot},
we get $G\succ H$ and $m(G,x)>m(H,x)$ for $x>t(G)$. Next, we suppose
that $G$ has $k$ connected components $G_1,\ldots,G_k$ where $k\geq
2$. Since $H$ is an edge-deleted subgraph of $G$, let $H$ be a union
of $k$ disjoint graphs $\{H_1,\ldots,H_k\}$, where $H_i$ is the
spanning subgraph of $G_i$. By Lemma \ref{union}, we know
\begin{equation}\label{1equ}
m(G,x)=m(\bigcup_{i=1}^kG_i,x)=\prod_{i=1}^km(G_i,x),
\end{equation}
\begin{equation}\label{2equ}
m(H,x)=m(\bigcup_{i=1}^kH_i,x)=\prod_{i=1}^km(H_i,x).
\end{equation}
Since $H$ is the proper spanning subgraph of $G$, we can find that
for some $i$, $H_i$ is the proper spanning subgraph of $G_i$.
Without loss of generality, assume that $H_1$ is the proper spanning
subgraph of $G_1$. By Proposition \ref{subroot}, we have $G_1\succ
H_1$. Then, for $x\geq t(G_1)$, $m(H_1,x)>m(G_1,x)\geq 0$. For other
subgraphs $H_i$ where $i\geq 2$, by Proposition \ref{Psub}, we have
$G_i\succeq H_i$. Then, for $x\geq t(G_i)$, $m(H_1,x)\geq
m(G_1,x)\geq 0$. From Lemma \ref{sub1}, we know that $t(G)\geq
max\{t(G_1,\ldots,t(G_k)\}$. Thus, for $x>t(G)$, we get
$m(H_1,x)>m(G_1,x)>0$ and $m(H_i,x)\geq m(G_i,x)>0$ for $2\leq i\leq
k$. Using Equations \ref{1equ} and \ref{2equ}, we deduce that for
$x>t(G)$,
\begin{equation*}
m(H,x)=\prod_{i=1}^km(H_i,x)>\prod_{i=1}^km(G_i,x)=m(G,x).
\end{equation*}

The proof is now complete. \qed

\begin{thm}\label{THM1}
Let $G$ be a connected graph. Assume that $G'$ is obtained from $G$
by a number of Kelmans transformations, and $G'$ is not isomorphic
to $G$. Then $G'\succ G$.
\end{thm}
\pf Let $G$ be a connected graph. Let $G_1$ be a graph obtained from
$G$ by the Kelmans transformation between $u$ and $v$, where $u$ is
the beneficiary, and $G_1$ is not isomorphic to $G$. We just need to
prove $G_1\succ G$.

Firstly, we find that $u$ has an adjacent vertex $w$ which is not
adjacent to $v$, and $v$ has an adjacent vertex $w'$ which is not in
the set of neighbors of the vertex $u$. Otherwise, $G_1$ is
isomorphic to $G$. This is a contradiction.

Secondly, we use Lemma \ref{deledge} and get
\begin{eqnarray*}
m(G,x)&=&m(G-uw,x)-m(G-u-w,x),\\
m(G_1,x)&=&m(G_1-uw,x)-m(G_1-u-w,x).
\end{eqnarray*}
Then we have the following equation
\begin{equation}\label{equation}
m(G,x)-m(G_1,x)=m(G-uw,x)-m(G_1-uw,x)+m(G_1-u-w,x)-m(G-u-w,x).
\end{equation}
Next, we consider the right part of the above equation. If we apply
the Kelmans transformation between $u$ and $v$ in $G-uw$, then we
can get the graph $G_1-uw$. Thus, by Theorem \ref{PCT}, we get
$G_1-uw\succeq G-uw$. That is,
\begin{equation*}
m(G-uw,x)-m(G_1-uw,x)\geq 0,
\end{equation*}
for $x\geq t(G_1-uw)$. Since $v$ has an adjacent vertex $w'$ and
$uw'$ is not in $E(G)$, $G_1-u-w$ is a proper spanning subgraph of
$G-u-w$. By Proposition \ref{Proper}, we have
\begin{equation*}
m(G_1-u-w,x)-m(G-u-w,x)>0,
\end{equation*}
for all $x>t(G-u-w)$. By Lemma \ref{sub1}, we find that $t(G_1)\geq
t(G_1-uw)$ and $G-u-w$ is isomorphic to a proper subgraph of $G_1$
by mapping $v$ of $G-u-w$ to $u$ of $G_1$. If $G_1$ is connected,
then we have $t(G_1)>t(G-u-w)$ by Lemma \ref{sub1}. Otherwise,
suppose that $G$ is connected and $G_1$ is not connected. Then it
follows that the set of neighbors of $u$ and the set of neighbors of
$v$ are disjoint. Moreover, $G_1$ is the disjoint union of a
isolated vertex $v$ and a connected component $H$. Then, it implies
that $G-u-w$ is isomorphic to a proper subgraph of $H$. Thus we have
$t(G_1)=t(H)>t(G-u-w)$ from Lemma \ref{sub1}. We can conclude that
$m(G-uw)-m(G_1-uw)\geq 0$ and $m(G_1-u-w,x)-m(G-u-w,x)>0$ for $x\geq
t(G_1)$. Finally, by Equation \ref{equation}, we have
$m(G,x)-m(G_1,x)>0$ that is $G_1\succ G$.

Now, it is time to prove the theorem. Suppose that $G'$ is obtained
from $G$ by $k$ Kelmans transformations
$\{KT(G,u,v),KT(G_1,u_1,v_1),\ldots,KT(G_{k-1},u_{k-1},v_{k-1})\}$
in order. From the above result, $G_1\succ G$. Moreover, from
Theorem \ref{PCT}, $G_2\succeq G_3,\cdots, G'\succeq G_{k-1}$. By
Proposition \ref{trans2}, we get $G'\succ G$. The proof is thus
complete. \qed

\section{The maximum matching roots of odd-cycle graphs}\label{sec4}

In this section, we first prove a theorem about the maximum matching
root of $F(n,m)$.
\begin{thm}\label{THM}
Given a positive integer $n$, for $n-1\leq m\leq \lfloor
\frac{3(n-1)}{2}\rfloor-1$, we have
$$t(F(n,m))<t(F(n,m+1)).$$
Given a positive integer $m$ with $m\geq 4$, for
$\lceil\frac{2m+1}{3}\rceil\leq n\leq m$, we have
$$t(F(n,m))<t(F(n+1,m)).$$
\end{thm}
\pf For a positive integer $n$, when $n-1\leq m\leq \lfloor
\frac{3(n-1)}{2}\rfloor-1$, $F(n,m)$ is a proper spanning subgraph
of $F(n,m+1)$. Then, by Proposition \ref{subroot}, we have
$t(F(n,m))<t(F(n,m+1)).$

Given a positive integer $m$ with $m>4$, when
$\lceil\frac{2m+1}{3}\rceil\leq n\leq m$, suppose that $G'$ is the
disjoint union of $F(n,m)$ and an isolate vertex named $w$. Since
$n\leq m$, it follows that there exists a cycle in $F(n,m)$. Due to
the definition of $F(n,m)$, the cycle in $F(n,m)$ is a triangle.
Then there exists a triangle $\{v_1,v_2,v_3\}$ where $v_1$ has the
maximum degree in $F(n,m)$. Up to isomorphism, the graph $F(n+1,m)$
can be obtained from $G'$ by deleting the edge $v_2v_3$ and adding
the edge $v_1w$. By Lemma \ref{deledge}, we have
\begin{eqnarray}
m(G',x)&=&m(G'-v_2v_3,x)-m(G'-v_2-v_3,x),\\
m(F(n+1,m),x)&=&m(F(n+1,m)-v_1w,x)-m(F(n+1,m)-v_1-w,x).
\end{eqnarray}

Note that The graph $G'-v_2v_3$ is isomorphic to $F(n+1,m)-v_1w$.
Then we get
\begin{equation}
m(G',x)-m(F(n+1,m),x)=m(F(n+1,m)-v_1-w,x)-m(G'-v_2-v_3,x).
\end{equation}

Thanks to the structure of $F(n+1,m)$, the graph $F(n+1,m)-v_1-w$ is
the disjoint union of $m-n$ edges and $3n-2m-1$ isolated vertices.
The graph $G'-v_2-v_3$ is the disjoint union of the graph
$F(n-2,m-3)$ and an isolated vertex $w$. By the definition of
$F(n-2,m-3)$, we know that the graph $F(n+1,m)-v_1-w$ is a subgraph
of $G'-v_2-v_3$. When $m>n$, $F(n+1,m)-v_1-w$ is isomorphic to a
proper spanning subgraph of $G'-v_2-v_3$. When $m=n$, since $m\geq
4$, we deduce that $F(n+1,m)-v_1-w$ is the disjoint union of $m-1$
isolated vertices, and the graph $G'-v_2-v_3$ is the disjoint union
of the graph $F(m-2,m-3)$ and an isolated vertex $w$. Thus
$F(n+1,m)-v_1-w$ is isomorphic to a proper spanning subgraph of
$G'-v_2-v_3$.

Above all, by Lemma \ref{Proper}, for $x>t(G'-v_2-v_3)$, we have
$m(F(n+1,m)-v_1-w,x)>m(G'-v_2-v_3,x)$. Since $G'-v_2-v_3$ is a
proper subgraph of $F(n+1,m)$, by Lemma \ref{sub1}, we get
$t(G'-v_2-v_3)<t(F(n+1,m))$. Then for $x\geq t(F(n+1,m))$, we obtain
$m(G',x)-m(F(n+1,m),x)>0$, that is, $F(n+1,m)\succ G'$. Thus, we
conclude that $t(F(n,m))=t(G')<t(F(n+1,m))$. The proof is now
complete. \qed

Next, we give a proof of Conjecture \ref{conjecture2} by proving the
following theorem.
\begin{thm}\label{THM2}
If $G$ is an odd-cycle graph of order $n$, then $t(G)\leq t(H_n)$,
and equality holds if and only if $G\cong H_n$.
\end{thm}
\pf If $G$ is disconnected, then there is a connected odd-cycle
graph $G'$ with order $n$ which contains $G$ as a proper subgraph.
By Lemma \ref{sub1}, $t(G)$ is less than $t(G')$. Then, suppose that
$G$ is a connected odd-cycle graph with order $n$ and size $m$.

By Kelmans transformations, we transfer the graph $G$ into $F(n,m)$.
By Theorem \ref{THM1}, we deduce that $t(G)\leq F(n,m)$ and equality
holds if and only if $G\cong F(n,m)$. When
$t<\lfloor\frac{3(n-1)}{2}\rfloor$, by Theorem \ref{THM} it follows
that $t(F(n,t))<t(F(n,\lfloor\frac{3(n-1)}{2}\rfloor))=t(H_n)$.
Then, we conclude that $t(G)\leq t(H_n)$ and equality holds if and
only if $G\cong H_n$. The proof is thus complete. \qed

Now, it is time to prove Theorem \ref{THM3}.

\textbf{Proof of Theorem \ref{THM3}.} For $m\leq 2$, it is easy to
check that the first two parts of the theorem is true.

For $m=3$, $G$ is one of the following graphs: the disjoint union of
three edges and $n-6$ isolated vertices, the disjoint union of the
star $K_{1,2}$, an edge and $n-5$ isolated vertices, the disjoint
union of the star $K_{1,3}$ and $n-4$ isolated vertices, the
disjoint union of a triangle and $n-4$ isolated vertices. By some
computations, it is easy to check that the third part of the theorem
is true.

For $4\leq m\leq n-2$, suppose that $G$ has $s$ connected components
$\{G_1,\cdots,G_s\}$ with order $\{n_1,\cdots,n_s\}$ and size
$\{m_1,\cdots,m_s\}$, where $s\geq n-m$. Since every connected
component is an odd-cycle graph, by the Kelmans transformations, we
can obtained a graph $G'$ with $s$ connected components
$\{F(n_1,m_1),\cdots,F(n_s,m_s)\}$. Then by Lemma \ref{union} and
Lemma \ref{THM1}, we have
$$t(G)=max\{t(G_1),\ldots,t(G_s)\}\leq max\{t(F(n_1,m_1)),\cdots,t(F(n_s,m_s))\}=t(G').$$

In every nontrivial connected component $F(n_i,m_i)$, one of the
vertices which have the maximum degree in $F(n_i,m_i)$ is named the
root vertex of $F(n_i,m_i)$. By identifying all root vertices, we
obtain a graph $F(n-s+1,m)$. Let $G''$ be the graph obtained from
the disjoint union of the graph $F(n-s+1,m)$ and $s-1$ isolated
vertices.

If there are at least two nontrivial connected components, then
every nontrivial connected component $F(n_i,m_i)$ is a proper
subgraph of $F(n-s+1,m)$. Then by Lemma \ref{union} and Lemma
\ref{sub1}, it implies that $t(G')<t(F(n-s+1,m))$. It follows that
$$t(G)\leq t(G')<t(G'').$$

Next, by Theorem \ref{THM}, for $s>n-m$ we have
$$t(F(n-s+1,m))<t(F(m+1,m)).$$
Above all, together with Theorem \ref{THM1}, we conclude that
$t(G)\leq t(F(m+1,m))$ and the equality holds if and only if $G$ has
only one nontrivial connected component which is isomorphic to
$F(m+1,m)$. Thus, the fourth part of the theorem is true.

When $n-1\leq m\leq \lfloor \frac{3(n-1)}{2}\rfloor$ and $m\geq 4$,
suppose that $G$ is disconnected and has $s$ connected components.
Then by the similar method to the fourth part, we can deduce that
the maximum matching root of $G$ is smaller than the maximum
matching root of $F(n-s+1,m)$. By Theorem \ref{THM}, for $s>1$ we
have
$$t(F(n-s+1,m))<t(F(n,m)).$$
Then, we conclude that $t(G)\leq t(F(n,m))$, and the equality holds
if and only if $G$ is connected and isomorphic to $F(n,m)$ from
Theorem \ref{THM1}. Thus, the fifth part of the theorem is true.

The proof is thus complete. \qed


\begin{thebibliography}{20}

\bibitem{ABS}A. Anuradha, R. Balakrishnan, W. So, Skew spectra of graphs
without even cycles, {\it Linear Algebra Appl.} 444(2014), 67--80.

\bibitem{BM}J.A. Bondy, U.S.R. Murty, {\it Graph Theory with Applications},
MacMillan Press Ltd., London and Basingstoke, 1978.

\bibitem{BH}A.E. Brouwer, W.H. Haemers, {\it Spectra of Graphs},
Springer, Berlin, 2012.

\bibitem{BCD}J. Brown, C. Colbourn, J. Devitt, Network transformation
and bounding network reliability, {\it Networks} 23(1993), 1--17.

\bibitem{CCF}M. Cavers, S.M. Cioab\v{a}, S.Fallat, D.A. Gregory,
W.H. Haemers, S.J. Kirkland, J.J. McDonald, M. Tsatsomeros,
            Skew-adjacency matrices of graphs,
            {\it Linear Algebra Appl.} 436(2012), 4512--1829.

\bibitem{CLL}L. Chen, X. Li, H. Lian, Lower bounds of the skew spectral
radii and skew energy of oriented graphs, arxiv:1405.4972v3.

\bibitem{PC1}P. Csikv\'{a}ri, On a conjecture of V. Nikiforov,
{\it Disc. Math.} 309(2009), 4522--4526.

\bibitem{PC}P. Csikv\'{a}ri, Applications of the Kelmans transformation:
extremality of the threshold graphs, {\it Electron. J. Combin.} 18(2011), \#P182.

\bibitem{FA}E.J. Farrell, An introduction to matching polynomials,
{\it J. Comb. Theory Ser. B}, 27(1979), 75--86.

\bibitem{GOD}C. Godsil, {\it Algebraic Combinatorics}, Chapman and Hall, New York 1993.

\bibitem{GR}C. Godsil, G. Royle, {\it Algebraic Graph Theory}, Graduate Texts in Mathematics,
Vol. 207, Springer--Verlag, New York, 2001.

\bibitem{GG}C. Godsil, I. Gutman, On the theory of the matching polynomial,
{\it J. Graph Theory} 5(1981), 137--144.

\bibitem{GUT}I. Gutman, A note on analogies between the characteristic and the
matching polynomial of a graph, {\it Publications de l'Institut
Math\'{e}matique} 31(1982), 27--31.

\bibitem{GMT}I. Gutman, M. Milun, N. Trinajsti\'{c}, Non-parametric resonance
energies of arbitrary conjugated systems. {\it J. Amer. Chem. Soc.} 99(1977), 1692--1704.

\bibitem{HL}O.J. Heilman, E.H. Lieb, Theory of monomer-dimer systems,
{\it Commun. Math. Physics} 25(1972), 190--232.

\bibitem{KEL}A.K. Kelmans, On graphs with randomly deleted edges,
{\it Acta. Math. Acad. Sci. Hung.} 37(1981), 77--88.

\bibitem{SSS}A. Satyanarayana, L. Schoppman, C.L. Suffel,
A reliability--improving graph transformation with applications
to network reliability, {\it Networks} 22(1992), 209--216.

\bibitem{SS}B. Shader, W. So, Skew spectra of oriented graphs,
{\it Electron. J. Combin.} 16(2009), \#N32.

\end{thebibliography}
\end{document}